\documentclass[11pt]{article}
\usepackage{amsmath}
\usepackage{latexsym,boxedminipage}
\newcommand{\Om}{\begin{array}{c}{\Omega}\vspace{-.1cm}\\{\ge}\end{array}}

\newcommand{\ld}{\left(}
\newcommand{\rd}{\right)}

\newtheorem{lemma}{Lemma}
\newtheorem{theorem}{Theorem}
\newtheorem{corollary}{Corollary}                                 
\newtheorem{proposition}{Proposition}                                 
\title{Partitions and Compositions Defined by Inequalities }
\author{Sylvie Corteel\\
{CNRS PRiSM, UVSQ}\\
 {45 Avenue des Etats-Unis}\\
{ 78035 Versailles, France}\\
{\tt syl@prism.uvsq.fr}\\
\and
 Carla D. Savage
\thanks{Research supported by NSA grants MDA 904-00-1-0059 and 
MDA 904-01-0-0083}\\
{Dept. of Computer Science}\\
 {N. C. State University, Box 8206 }\\
{ Raleigh, NC 27695, USA}\\
{\tt savage@csc.ncsu.edu}                                                    
}

\begin{document}
\maketitle
\begin{abstract}
We consider sequences of integers
$(\lambda_1, \ldots, \lambda_k)$ defined by a system of linear 
inequalities 
$\lambda_i \geq \sum_{j>i} a_{ij} \lambda_{j}$
with integer  coefficients.
We show that when the constraints are strong enough to guarantee that all
$\lambda_i$ are nonnegative, the generating function for the integer
solutions of weight $n$
has a finite product
form $\prod_{i} (1-q^{b_i})^{-1}$, where the $b_i$ are
positive integers that can be computed from the
coefficients of the inequalities.
The results are proved bijectively and are used to give several
examples of interesting identities for integer partitions and
compositions.
The method can be adapted to accommodate equalities along with
inequalities and can be used to obtain multivariate forms of
the generating function.
We show how to extend the technique to obtain the generating function
when the 
coefficients $a_{i,i+1}$ are allowed to be rational,
generalizing the case of lecture hall
partitions.
Our initial results were
conjectured thanks to the Omega package \cite{Om3}.
\end{abstract}

\section{Introduction}
For a sequence $\lambda = (\lambda_1, \lambda_2, \ldots, \lambda_k)$
of integers,
define the {\em weight} of $\lambda$ to be 
$\lambda_1 + \cdots + \lambda_k$ and call each $\lambda_i$
a {\em part} of $\lambda$.  If a sequence $\lambda$ of weight $n$
has all parts nonnegative, we call it a {\em composition} of $n$
into $k$ nonnegative parts and if, in addition, $\lambda$ is a 
nonincreasing sequence, we call
it a {\em partition} of $n$ into at most $k$  parts. In the sequel
we will consider that $\lambda_i=0$ if $i<0$ or $i>k$.

In this paper we  study partitions and
compositions into $k$ nonnegative 
parts defined by equalities and inequalities. 
This work was motivated by
results of the form ~:
\begin{itemize} 
\item Given a positive integer $r$, the
partitions $\lambda= \lambda_1, \ldots, \lambda_k$ of $n$
which satisfy
$\lambda_i \geq r \lambda_{i+1}$ for $1 \leq i \leq k$ 
have weight
generating function  
$\prod_{i=0}^{k-1} \left(1-q^{1+r+ \cdots + r^i}\right)^{-1}$ 
\cite{hic}.

\item Given a positive integer $r$ the 
generating function of the partitions 
$\lambda$ of $n$
with at most $k$ parts and  $\lambda_i\ge \sum_{j=1}^{r}
(-1)^{j+1}\binom{r}{j}\lambda_{j+i}$, $1\le i<k$ is~:
$\prod_{i=0}^{k-1} \left(1-q^{\binom{i+r}{r}}\right)^{-1}.$ 
See \cite{Om2,CCH,SC}.

\item The generating function of the partitions 
$\lambda$ of $n$ with at most $k$ parts and  
$\lambda_i/(k-i+1)\ge \lambda_{i+1}/(k-i)$, 
$1\le i<k$ is~:
$\prod_{i=1}^{k} \left(1-q^{2i-1}\right)^{-1}.$ See the
{\em Lecture Hall Theorem } in \cite{BME1}. 
\end{itemize}
More generally, we consider integer sequences 
$\lambda$ of length $k$
satisfying $\lambda_i\ge \sum_{j=i+1}^{k}
a_{ij}\lambda_{i}$ where the $a_{ij}$ guarantee that
all $\lambda_i\ge 0$.
We show in Section 2 that when the $a_{ij}$ are all integers,
the generating function for these compositions   
is $\prod_{i=1}^{k} (1-q^{b_i})^{-1}$, where  
$b=(b_1,\ldots ,b_{k})$ is
a sequence of positive integers.  Several generalizations are included.

In Section 3, we consider rational coefficients $a_{ij}$.
We show how to use the results of Section 2 to give an explicit form
for the generating function for {\em any} set of compositions defined by
the ratio of consecutive parts:
\begin{equation*}
\lambda_1  \geq  \frac{n_1}{d_1}\lambda_2;\  \ \ \
\lambda_2  \geq   \frac{n_2}{d_2}\lambda_3;\  \  \ \
\lambda_3  \geq   \frac{n_3}{d_3}\lambda_4;\  \  \ \
 \ldots;  \ \ \
\lambda_{k-1}  \geq   \frac{n_{k-1}}{d_{k-1}}\lambda_k;\  \   \ \
\lambda_{k}  \geq  0.
\end{equation*}
This result has been implemented in Maple and our experiments have led
to several interesting results.  
In Section 3.2 we focus on one special case which has consequences such
as:
the compositions of $n$ into at most $k$ parts satisfying
\[
\lambda_1 \geq \lambda_2 \geq \lambda_3/2 \geq \lambda_4 \geq \lambda_5/2
\geq \cdots \geq \lambda_{2k} \geq \lambda_{2k+1}/2 \geq 0
\]
are in one-to-one correspondence with the set of partitions of $n$
into parts of size at most $3k+2$ in which parts divisible by 3 can
appear at most once.


In \cite{BME1}, Bousquet-M{\'e}lou and Eriksson considered the set $L_k$
of partitions $\lambda$ into at most $k$ parts satisfying
$\lambda_i/(k-i+1)\ge \lambda_{i+1}/(k-i)$, for $1\le i<k$,
and proved the following {\em Lecture Hall Theorem}:
\begin{equation}
\sum_{\lambda \in L_k} q^{|\lambda|} = 
\prod_{i=1}^{k}\frac{1}{1-q^{2i-1}}.
\label{bmelhp1}
\end{equation}
In Section 4, we study the 2-variable generating function
and show a slight generalization of (\ref{bmelhp1}) in 
which the constraints on $\lambda_1$ can be modified.

Our initial results were conjectured thanks to experiments with the Omega
 package
\cite{Om3,Om6}, a
Mathematica  implementation
of the Omega operator defined by MacMahon \cite{MaMo}.
This operator was then not used for 85 years except by Stanley 
in 1973 \cite{St73}. A few years ago Andrews
revived this operator \cite{Om1,Om2} and used it in
\cite{Om1} to give a second
proof of the Lecture Hall Theorem.
In conjunction with Paule and Riese, he implemented the operator in the
Omega package and together they have  continued to identify
the power of the Omega operator for such
combinatorial problems as 
 magic squares \cite{Om5},
hypergeometric multisums \cite{Om4}, 
constrained compositions \cite{Om7}, plane partitions diamonds 
\cite{Om8},
and 
$k$-gons partitions \cite{Om9}.

\section{Integer Coefficients}

\subsection{Overview}

A matrix $A$ is called  {\em upper triangular} if $A[i,j]=0$ for
$i>j$ and {\em strictly upper triangular} if  $A[i,j]=0$ for
$i \geq j$.

Let $A[1..k,1..k]$ be  a strictly upper triangular matrix  of
integers  and
let $P_A$ be the set of sequences
$\lambda = (\lambda_1, \lambda_2, \ldots, \lambda_k)$
satisfying
\begin{equation}
\lambda_i  \geq   \sum_{j=i+1}^{k} A[i,j] \lambda_{j} \ \ \
{\mbox{\rm  for $1 \leq i \leq k$}}.
\label{inequalities}
\end{equation}

In this section
we show that if every $\lambda \in P_A$ is a sequence of nonnegative integers,
 that is, a composition,
then the generating function for the elements of $P_A$ of weight $n$
is
\[
\prod_{i=1}^{k}\frac{1}{1-q^{b_i}}
\]
where $\{b_i\}$ is  sequence of positive integers
which can be computed from  $A$.

\subsection{The Details}

Given
$A[1..k,1..k]$,  a strictly  upper triangular matrix of
integers,
we know  $A$ is nilpotent, that is,  $A^k=0$.
 Considering $\lambda$ as a column vector, the matrix inequality,
$\lambda \ge A \lambda,$ describes the solutions to (\ref{inequalities}).
 Then
\begin{equation}
\lambda =  A\lambda + [s_1, s_2, \ldots s_{k}]^T
\label{slack}
\end{equation}
where $s =(s_1, \ldots , s_k)$ is the  column vector of nonnegative integers
defined by
\[
s_i = \lambda_i - \sum_{j=i+1}^k A[i,j] \lambda_j.
\]
Iterating (\ref{slack})  yields
\[
\lambda =(I+A+A^2+\cdots +
A^{k-1})  [s_1, s_2, \ldots s_{k}]^T \]
since $A^k=0$.
Since $(I-A)(I+A+A^2 + \cdots + A^{k-1})=I$, 
$I-A$ is invertible, 
and its  inverse is an upper triangular
matrix of integers, since $A$ is.
Thus,
\[
\lambda =  (I-A)^{-1} [s_1, s_2, \ldots s_{k}]^T.
\]
Let $B=(I-A)^{-1}$. Then
\begin{equation}
P_A = \{B[s_1, \ldots, s_k]^T \ | \ s_i \geq 0  \ {\rm for} \ 
1 \leq i \leq k\}.
\label{PAdef}
\end{equation}

\begin{lemma}
The set $P_{A}$, of solutions to
{\rm (\ref{inequalities})}, is a set of compositions
if and only if every element
of $B=(I-A)^{-1}$ is nonnegative.
\label{PAlemma}
\end{lemma}

\noindent
{\bf Proof}
If $\lambda \in P_A$, then from (\ref{PAdef}),
$\lambda=B[s_1,\ldots, s_k]^T$ for some sequence of nonnegative integers
$\{s_i\}$. If all elements of $B$ are nonnegative, $\lambda_i \geq 0$
for $1 \leq i \leq k$.  However, if some entry $B[j,l]$ is negative,
then choosing $\{s_i\}$ so that $s_l = 1$ and $s_i=0$ if $i \not = l$
makes entry $\lambda_j$ of $\lambda=B[s_1,\ldots, s_k]^T$ negative.
\hfill $\Box$

\begin{theorem}
Let
$P_{A}(n,k)$
 be the set of sequences
$\lambda = (\lambda_1, \lambda_2, \ldots, \lambda_k)$
of weight $n$  satisfying
\begin{equation*}
\lambda_i  \geq   \sum_{j=i+1}^{k} A[i,j] \lambda_{j} \ \ \
{\mbox{\rm  for $1 \leq i \leq k$}}.
\end{equation*}
If every $\lambda \in P_{A}$ is a composition,
then
the generating function for $P_{A}(n,k)$ is
\begin{eqnarray}
\sum_{n=0}^{\infty}|P_{A}(n,k)|q^n \ = 
\prod_{i=1}^{k}\frac{1}{1-q^{b_i}},
\label{gf=}
\end{eqnarray}
where $b_i = \sum_jB[j,i]$ and $B=(I-A)^{-1}$.
\label{PAtheorem}
\end{theorem}

\noindent
{\bf Proof}
Let $h=[1, 1, \ldots , 1] $ be the row vector of length $k$
containing only ones. The weight of a composition 
$\lambda \in P_A$
is
$|\lambda| = h\lambda =h B  [s_1, \ldots s_{k}]^T$,
using (\ref{PAdef}),
where
$s_i = \lambda_i - \sum_{j=i+1}^k A[i,j] \lambda_j$.
Define $b=hB$. Then
 $b_i=\sum_j B[j,i]$ and
 the
weight generating function of $P_A$ is
\begin{equation}
\sum_{\lambda \in P_A}q^{\lambda_1+ \cdots + \lambda_k} =
\sum_{s_1, \ldots , s_{k} \ge 0} q^{b_1s_1+\cdots + b_ks_{k} }
=\prod_{i=1}^k \frac 1{1-q^{b_i}}.
\label{linalg}
\end{equation}
\hfill $\Box$

Note that this argument establishes that the mapping
\[
\Theta(\lambda) = (b_{1}{s_{1}}, \ldots, b_{k}{s_{k}})
\]
where
\[
s_i = \lambda_{i}  - \sum_{j=i+1}^{k}A[i,j] \lambda_{j}
\]
is a bijection from $P_{A}(n,k)$ to the set of sequences
$r_1, \ldots, r_{k}$ of weight $n$ in which $r_i$ is a nonnegative
multiple of
$b_i$.  In the case that the $b_i$ are distinct positive integers, this
can be viewed alternatively as a bijection with partitions of $n$ into
parts in $\{b_1, \ldots, b_{k}\}$.

\noindent {\bf Example 0}
Starting with  the matrix
$$A=\left[\begin{array}{cccc}
0& 2 & -1 & 0 \\
0 & 0 & 2 & -1\\
0 & 0 & 0 & 2\\  
0 & 0 & 0 & 0\\  
\end{array}\right]\ \ \  {\rm we}\ {\rm get}\ \ \  
B=(I-A)^{-1}=\left[\begin{array}{llll}
1 & 2 & 3 & 4\\
0 & 1 & 2 & 3\\
0 & 0 & 1 & 2 \\ 
0 & 0 & 0 & 1 \\ 
\end{array}\right]$$
and summing each column of $B$ gives the sequence $b=(1,3,6,10)$.
By Theorem \ref{PAtheorem}, the generating function  
for $P_{A}$ is $[(1-q)(1-q^3)(1-q^6)(1-q^{10})]^{-1}$.

If $A$ is a strictly upper triangular matrix of integers
with the property that
$P_A$ is a set of partitions and if the $b_i$ are distinct
positive integers, Theorem 1 gives a bijection between two sets of
partitions.   Several such sample results, including  Example 1 below,
appear throughout this section.  In other situations, such as Example 2
below,
Theorem 1 gives a bijection between families of compositions.
Later examples, such as Example 15 in Section 3,
 will show bijections between a family of compositions,
on the one hand, and a family of partitions on the other.

\noindent
{\bf Example 1}
The partitions $\lambda$ of $n$
with at most $k$ parts and with $\lambda_1 \geq \sum_{i=2}^{k}
\lambda_{i}$  comprise $P_A$ where, 
$A[i,j]=1$ if $1=i<j$ or $j=i+1$ and $A[i,j]=0$ otherwise.
Then the matrix $B=(I-A)^{-1}$ is:
\[
B[i,j]=\left\{
\begin{array}{ll}
0 & {\rm if}\ j<i\\
1  & {\rm if}\ j=i \ {\rm or}\ j>i>1\\
j-1 & {\rm otherwise.}
\end{array}\right.
\]
Summing the $i$th column in the matrix $B=(I-A)^{-1}$ gives
$b_i=2(i-1)$ for $2 \leq i \leq k$ and $b_1=1$, so by Theorem 1,
the generating function is
$$(1-q)^{-1}\prod_{i=1}^{k-1} (1-q^{2i})^{-1}.$$



Let $P_{A,=}(n,k)$ be the set of compositions of $n$ in $P_A$ 
with at most $k$ parts and with 
$\lambda_1  =  \sum_{j=2}^{k} A[i,j] \lambda_{j}$.
Since this forces $s_0=0$ in (\ref{PAdef}), we get the following
from \ref{linalg}.

\begin{corollary}
Let $P_{A,=}(n,k)$ be the set of compositions in $P_A$ of $n$
with at most $k$ parts and with 
$\lambda_1  =  \sum_{j=2}^{k} A[i,j] \lambda_{j}$.
The generating function for $P_{A,=}(n,k)$ is
\[
\prod_{i=2}^{k}\frac{1}{1-q^{b_i}}.
\]
\label{PAeq}
\end{corollary}


\noindent
{\bf Example 2}
There is a one-to-one correspondence between sequences
$\lambda_1, \ldots, \lambda_k$ of weight $n$ satisfying
\begin{eqnarray*}
\sum_{j=0}^{k-1}(-1)^{j}\lambda_{1+j} &  = &  0\\
\sum_{j=0}^{k-i}(-1)^{j}\lambda_{i+j} &  \geq & 0   \ \ \ \ \ \ 
\mbox{for $2 \leq i \leq k$}
\end{eqnarray*}
and the set of compositions of $n$ into $k-1$ even parts.
The $\lambda$ satisfying the constraints are compositions
and they comprise the set
$P_A(n,k)$ where $A$ is the $k \times k$
upper triangular matrix defined by
$A[i,i]=0$ and $A[i,j]=(-1)^{i+j+1}$.  The result follows from
Corollary \ref{PAeq} by summing the columns of $B=(I-A)^{-1}$
to get
$ b_2=b_3= \cdots = b_{k}=2$.
 
\noindent
{\bf Example 3}
 The generating function of the partitions $\lambda$ of $n$
with at most $k$ parts and with $\lambda_1\ge2\lambda_2+\sum_{i=2}^{k-1}
\lambda_{i+1}$ and $\lambda_2\ge\sum_{i=2}^{k-1}
\lambda_{i+1}$ is~:
$$\frac{1}{(1-q)(1-q^3)}\prod_{i=1}^{k-2} \frac{1}{1-q^{5i}}.$$
This follows from Theorem \ref{PAtheorem}:
Since $\lambda$ is a partition, for $i>2$, $A[i,j]=1$ when $j = i+1$ and
0 otherwise.
When $i \leq 2$,  $A[1,2]=2$ and otherwise $A[i,j]=1$ for
$i=1,2$, $j  >i$, so $b_1=1$, $b_2=3$, and $b_i=5(i-2)$ for
$3 \leq i \leq k$.

\subsection{Further Examples and Generalizations}

Theorem \ref{PAtheorem} gives bijective proofs of  identities
such as
the two first
examples given in the introduction.

\noindent
{\bf Example 4}  Given a positive integer $r$, the
partitions $\lambda= \lambda_1, \ldots, \lambda_k$ of $n$
which satisfy
$\lambda_i \geq r \lambda_{i+1}$ for $1 \leq i \leq k-1$ 
correspond to the constraint matrix $A$ with
$A[i,i+1]=r$.
Then $B$ has $B[i,j]=r^{j-i}$, $j \geq i$, and 0 otherwise,
 so by Theorem \ref{PAtheorem},
 the generating function is 
$\prod_{i=0}^{k-1} \left(1-q^{1+r+ \cdots + r^i}\right)^{-1}.$ 
See \cite{hic}. 


\noindent
{\bf Example 5}   Given $r$, a positive integer.
Consider the partitions $\lambda$ of $n$
with at most $k$ parts and 
 $\lambda_i\ge \sum_{j=1}^{r}
(-1)^{j+1}\binom{r}{j}\lambda_{j+i}$, that is~:
$A[i,j]=\left\{
\begin{array}{ll}
0 & {\rm if}\ i\le j\\
(-1)^{j-i+1}{r\choose j-i} &  {\rm otherwise}
\end{array}\right.$.
The generating function is
\[
\prod_{i=0}^{k-1}\frac 1{1-q^{r+i\choose r}},
\]
since it can be shown that 
$B[i,j]=\left\{
\begin{array}{ll}
0 & {\rm if}\ i  >j\\
{r-1+j-i\choose r-1} &  {\rm otherwise}
\end{array}\right. .$
See \cite{Om2,CCH}.

\noindent
{\bf Example 6} 
  The generating function of the partitions $\lambda$ of $n$
with at most $k$ parts and  $\lambda_i\ge\lambda_{i+1}
+\lambda_{i+2}$, $1\le i\le k$ is~:
$\prod_{i=2}^{k+1} (1-q^{F_{i}-1})^{-1}$, where $F_i$ is the
 $i^{th}$ Fibonacci number
(defined by  $F_0=F_1=1$ and $F_i=F_{i-1}+F_{i-2}$
for $i>1$.)
In this case, $B[i,j]=F_{j-i}$ if $j \geq i$ and 0 otherwise.


\noindent
{\bf Example 7}
Studying partitions with  $\lambda_i\ge i \sum_{j=i+1}^{k}
\lambda_{j}$, $1\le i\le k$, we get
$B[i,j]=1$ if $i=j$,  $j!/((i+1)(i-1)!)$ if $i<j$ and
0 otherwise. It is easy to show by induction that
\[
\sum_{i=1}^{j-1}\frac{1}{(i+1)(i-1)!}=\frac{j!-1}{j!}.
\]
Hence $b_j=\sum_{i}B[j,i]=1+\sum_{i=1}^{j-1}\frac{j!}{(i+1)(i-1)!}=j!$.
The generating function of the partitions $\lambda$ of $n$
with at most $k$ parts and $\lambda_i\ge i \sum_{j=i+1}^{k}
\lambda_{j}$, $1\le i\le k$, is~:
$\prod_{j=1}^{k} (1-q^{j!})^{-1}$.


We can generalize Theorem 1 to allow the constraints of the
matrix $A$ to be satisfied with equality for any  specified set
of $\lambda_i$.
                         
\begin{corollary}
Given a set $S \subseteq \{1,2, \ldots, k\}$, let
$P_A(n,k;S)$ be the set of  sequences
$\lambda = (\lambda_1, \lambda_2, \ldots, \lambda_k)$ 
of weight  $n$  satisfying, for $1 \leq i \leq k$:
\[
\lambda_i = \sum_{j=i+1}^{k} A[i,j] \lambda_{j} \ \ \ 
{\mbox{\rm  if $i \in S$}}
\]
\[
\lambda_i \geq \sum_{j=i+1}^{k} A[i,j] \lambda_{j} \ \ \ 
{\mbox{\rm  if $i \not\in S$}}.
\]
If all elements of $P_A$ are compositions,
the generating function for $P_A(n,k;S)$ is
\[
\sum_{n=0}^{\infty}|P_A(n,k;S)|q^n \ = 
\prod_{i=1, i \not\in S}^{k}\frac{1}{1-q^{b_i}},
\]
where the $b_i=\sum_j B[j,i]$.
\label{anyS}
\end{corollary}

\noindent
{\bf Proof}

For $1 \leq i \leq k$, let 
\[
s_i = \lambda_{i} - 
\sum_{j=i+1}^{k} A[i,j] \lambda_{j}.
\]
Any $\lambda \in P_A(n,k;S)$ must have 
$s_i=0$ for $i \in S$, so the
result follows from (\ref{linalg}) in the proof of Theorem \ref{PAtheorem}.
\hfill $\Box$

\noindent
{\bf Example 8}
Consider
 the set of  partitions $\lambda$ of $n$
with at most $k$ parts and with $\lambda_1=\sum_{i=2}^{k}
\lambda_{i}$ and $\lambda_2=\sum_{i=3}^{k}
\lambda_{i}$.
Here $S=\{1,2\}$, $b_1=1$, $b_2=2$, and
$b_i=4(i-2)$ for $3 \leq i \leq k$, so by
Corollary 2, the generating function is
$$ \prod_{i=1}^{k-2} (1-q^{4i})^{-1}.$$


\begin{corollary}
Given a set $S \subseteq \{1, \ldots k\}$
and nonnegative integers
$d_1, d_2, \ldots d_{k}$, let
$D_A(n,k;S)$ be the set of  sequences
$\lambda = (\lambda_1, \lambda_2, \ldots, \lambda_k)$ 
of weight $n$ into nonnegative parts satisfying, for $1 \leq i \leq k$:
\[
\lambda_i = \sum_{j=i+1}^{k} A[i,j] \lambda_{j} + d_{i} \ \ \ 
{\mbox{\rm  if $i \in S$}}
\]
\[
\lambda_i \geq \sum_{j=i+1}^{k} A[i,j] \lambda_{j} + d_{i} \ \ \ 
{\mbox{\rm  if $i \not\in S$}}.
\]
If all elements of $P_A$ are compositions,
the generating function for $D_A(n,k;S)$ is just
$$\prod_{i=1}^{k}q^{d_ib_i}
\sum_{n=0}^{\infty}|P_A(n,k;S)|q^n.$$
\label{add_d}
\end{corollary}

\noindent
{\bf Proof}
In this case,
(\ref{slack}) becomes
$\lambda \geq A \lambda + [s_1+d_1, \ldots, s_k+d_k]^T$,
where $s_i=0$ when $i \in S$.
By (\ref{linalg}) in the proof of Theorem \ref{PAtheorem},
the generating function becomes
\begin{equation*}
\sum_{\lambda \in D_A(n,k;S)}q^{\lambda_1+ \cdots + \lambda_k} =
\sum_{s_1, \ldots , s_{k} \ge 0} q^{b_1(s_1+d_1)+\cdots + b_k(s_{k}+d_k) }
=\prod_{i=1}^{k}q^{d_ib_i}
\sum_{n=0}^{\infty}|P_A(n,k;S)|q^n.
\end{equation*}
\hfill $\Box$

\noindent
{\bf Example 9}
Consider the partitions $\lambda_1, \ldots, \lambda_k$ of $n$
satisfying $\lambda_1=\lambda_2+\sum_{i=4}^{k}
\lambda_{i}$ or $\lambda_1=\lambda_2+\sum_{i=4}^{k}
\lambda_{i}+1$.  These are in one-to-one correspondence with
the partitions of $n$ into odd parts of size at most $2k-3$.
To see this, note that  for both sets of constraints, 
$S=\{1\}$,
$b_1=1, b_2=2$ and $b_i= 2i-3$
for $3 \leq i \leq k$.
For the first set, all $d_i=0$.
For the second set,  $d_1=1$, and $d_i=0$ for $i>1$.
  So by Corollary 3, the generating function is:

$$\frac{1+q}{(1-q^2)}\prod_{i=3}^{k-1}\frac{1}{ (1-q^{2i-3})} =
\prod_{i=1}^{k-1}\frac{1}{ (1-q^{2i-1})}.$$

\noindent
{\bf Example 10}  Partitions of $n$ into $k$ odd parts can be viewed as
those sequences $(\lambda_1, \ldots, \lambda_{2k})$ satisfying
$\lambda_{2k} \geq 0$ and for $i < 2k$,
\[
\lambda_i = \lambda_{i+1}+1 \ \ \ \  \mbox{if $i$ is odd}; \ \  \ \ 
\lambda_i \geq \lambda_{i+2} \ \ \ \  \mbox{if $i$ is even}.
\]
(Just let $\alpha_i=\lambda_{2i-1}+\lambda_{2i}$.)
For this system, $d_i = 1$ if $i$ is odd and $0$ if $i$ is even and
$S=\{1, 3, 5, \ldots, 2k-1\}$.
So, summing the columns of $B={I-A}^{-1}$,
 we get $b_i = 1$ if $i$ is odd and $i$ if $i$ is even and
by Corollary 3, the generating
function is
\[
q^k \prod_{i=1}^{k}(1-q^{2i})^{-1},
\]
giving (not surprisingly)
 a bijection with partitions of $n-k$ into $k$ nonnegative even parts.

\noindent{\bf Remark}
Note that we can link these results 
to partition analysis. The Omega operator $\Om$ \cite{MaMo}
is defined as follows~:
\[
\Om \sum_{s_1=-\infty}^{\infty}\ldots \sum_{s_r=-\infty}^{\infty}\
A_{s_1,\ldots ,s_r}\omega_1^{s_1}\ldots \omega_r^{s_r}=
\sum_{s_1=0}^{\infty}\ldots \sum_{s_r=0}^{\infty}\
A_{s_1,\ldots ,s_r}
\]
To calculate with this operator MacMahon \cite{MaMo} proposed a list of  
elimination rules. Here is one of them~:
\[
\Om \frac{1}{(1-\omega x)\ld 1-\frac{y}{\omega^r}\rd}=
\frac{1}{(1- x)(1-x^ry)},\ \ \ \ r\ge 0
\]
Our results can be translated in general elimination
rules. Let $k$ be a positive integer and $s=(s_1,\ldots,s_k)$
be a sequence of integers. Then ~:
\[
\Om \frac{1}{1-\omega x}
\prod_{i=1}^k \frac{1}{1-{y_i}/{\omega^{s_i}}}=
\frac{1}{1-x}\prod_{i=1}^k \frac{1}{1-x^{s_i}y_i},
\]
\[
\Om \frac{1}{1-x\omega_1/\omega_2}
\prod_{i=1}^k  \frac{1}{1-{y_i}/(\omega_1/\omega_2)^{s_i}}=
\prod_{i=1}^k \frac{1}{1-x^{s_i}y_i}.
\]

\subsection{From Product to Constraint Matrix}

For any sequence $c=(c_1, \ldots, c_k)$ of positive
integers,  can one construct a strictly upper triangular constraint matrix
$A[1..k,1..k]$ of integers such 
that $P_A$ is a set of compositions
with weight generating function
$\prod_{i=1}^{k}(1-q^{c_i})^{-1}$? We can answer yes for any sequence.
First let us suppose that $c_1=1$.

\begin{proposition}
Given a sequence  $(c_1, \ldots ,c_k)$ of positive integers, such that
$c_1=1$, 
there exist
\[
\prod_{i=2}^k {c_i-1\choose i-2}
\]
upper triangular matrices $B$ with nonnegative integers coefficients
and ones on the diagonal
such that if $A=I-B^{-1}$ then   $P_A$ has 
generating function $\prod_{i=1}^k(1-q^{c_i})^{-1}$.
\end{proposition}

\noindent
{\bf Proof }
Given a sequence $(c_1, \ldots ,c_k)$ of positive integers with $c_1=1$, 
we can always construct an
upper triangular matrix $B$ with ones on the diagonal and nonnegative
entries such that the sum of the entries in the $j^{th}$ column is
$c_j$.
Then $B$ is invertible and its inverse is an upper triangular, integer
matrix with ones on the diagonal.
Thus the matrix $A=I-B^{-1}$ is the constraint
matrix $A$ such that $P_A$ has
generating function $\prod_{i=1}^k(1-q^{c_i})^{-1}$.
\hfill{$\Box$}

\noindent{\bf Remark} If $c_1$ is not equal to 1, we can construct
matrices 
$A[1..k+1,1..k+1]$ of integers such 
that $P_{A,=}$ is a set of compositions
with weight generating function
$\prod_{i=1}^{k}(1-q^{c_i})^{-1}$.

When the $c_i$ satisfy $c_i \geq i$, we can use the following
corollary to show, in Proposition 2, a simple form for the
inequalities defining $P_A$.

\begin{corollary}

If $P_A$ is a set of compositions with $\lambda_i \geq \lambda_{i+1}$
for $2 \leq i \leq k-1$ and with  the first part
constrained by
$\lambda_1 \geq \sum_{i=1}^{k-1} a_i\lambda_{i+1}$,
the generating function of $P_A$ is
\[
\frac{1}{1-q}
\prod_{i=1}^{k-1} \frac{1}{1-q^{ i+a_1+ \cdots +a_i}}.
\]

\end{corollary}
\noindent
{\bf Proof}
In this case $B=(I-A)^{-1}$ is the upper triangular matrix defined by
$B[1,1]=1$,
$B[1,j]=a_1+ \dots +a_{j-1}$ for $j>1$, and for $2 \leq i \leq j$,
$B[i,j]=1$.
The column sums, $b_j$, satisfy
$b_1=1$ and for $2 \leq i \leq k$, $b_j=(j-1)+ a_1+ \cdots +a_{j-1}$,
so the result now follows from Theorem \ref{PAtheorem}.
\hfill{$\Box$}

\begin{proposition}
For any  sequence $c=(c_1, \ldots, c_k)$ of positive
integers satisfying $ c_i \geq {i-1}$, one can construct
a matrix $A$ such that either $P_A(n,k)$
or $P_{A,=}(n,k+1)$ is a set of
compositions with weight generating function
$\prod_{i=1}^{k}(1-q^{c_i})^{-1}$.
\end{proposition}

\noindent
{\bf Proof}
If $c_1=1$, let $A$ be the $k \times k$ strictly upper
triangular matrix such that
$A[i,i+1]=1$ if $i \geq 2$;
$A[i,j]=0$ if $2 \leq i < j-1$; 
$A[1,2]=c_2-1$; and  $A[1,j]=c_{j}-c_{j-1}-1$ for $j>2$.
Then by Corollary 4, using $a_i = A[1,i+1]$, $P_A(n,k)$ has generating function
$\prod_{i=1}^{k}(1-q^{b_i})^{-1}$,
where the $b_i$ are defined by
$b_1=1=c_1$ and for $2 \leq i \leq k$
\[
b_i= (i-1)+ \sum_{j=1}^{i-1}A[1,j+1]=
(i-1)+ (c_2-1) + \sum_{j=2}^{i-1}(c_{j+1}-c_{j}-1) = c_i.
\]
Note from the proof of Corollary 4 that
in this case $B=(I-A)^{-1}$ has all entries nonnegative,
since $\sum_{i=1}^{j-1}a_i=c_{i}-(i-1) \geq 0$, by the constraints on
the $c_i$.  So, by Lemma \ref{PAlemma}, $P_A(n,k)$ is a set of
compositions.

If $c_1>1$,  apply the technique for $c_1=1$ to the sequence
$c'=(1,c_1,c_2, \ldots c_k)$ to get a family of compositions
$P_A(n,k)$ with generating function
$(1-q)^{-1}
\prod_{i=1}^{k} (1-q^{c_i})^{-1}$.
Then $P_{A,=}(n,k)$ is the required set of compositions.
\hfill{$\Box$}

\noindent
{\bf Example 11}
Given the sequence  $1,3,5, \ldots 2k-1,$ of the first $k$ odd
positive integers,  the method
of Proposition 2 says that
$\prod_{i=1}^{k} (1-q^{2i-1})^{-1}$
 is the generating function of the partitions $\lambda$ of $n$
with at most $k$ parts and with $\lambda_1\ge 2\lambda_2+\sum_{i=3}^{k-1}
\lambda_{i}$.


\noindent
{\bf Example 12 }
Given the sequence $1,4,6,9,11,14,...,5k-4,5k-1$, of the first $2k$ integers
congruent to 1 or 4 mod 5, the method of Proposition 1 says that
if we choose $B$ such that
\[
B[i,j]=\left\{
\begin{array}{ll}
0 & {\rm if}\ j<i\\
1  & {\rm if}\ j=i\\
3  & {\rm if}\ i\ {\rm odd}\\
2 & {\rm otherwise.}
\end{array}\right .
\]
Then 
$$\prod_{i=1}^{k} \frac{1}{(1-q^{5i-4})(1-q^{5i-1)}}$$
is the 
  generating function of the partitions $\lambda$ of $n$
with at most $2k$ parts and with $\lambda_{2i-1}\ge\sum_{j\ge i}3(2^{j-i})
(\lambda_{2j}-\lambda_{2j+1})$ and $\lambda_{2i}\ge2\lambda_{2i+1}-
\sum_{j> i}2^{j-i+1}
(\lambda_{2j}-\lambda_{2j+1})$.

\noindent
{\bf Example 13}
Given the sequence $2,2,4,4,6,6,...,2k,2k$, by Proposition 2,
\[
\prod_{i=1}^{k}\frac{1}{(1-q^{2i})^2}
\]
is the generating function of the set of compositions satisfying
$\lambda_1 \geq \sum_{i=2}^{2k}(-1)^i \lambda_i$ and
$\lambda_i \geq \lambda_{i+1}$ for $2 \leq i \leq 2k$.

\section{Rational Coefficients}

In this section we would like to generalize our results to allow
some of the elements of the constraint matrix to be rational.
In particular, we will find the generating function for the set of
integer sequences $\lambda_1, \ldots, \lambda_k$ satisfying the
constraints:
\begin{eqnarray}
\lambda_1 & \geq & c_1 \lceil \frac{n_1}{d_1}\lambda_2 \rceil +
\sum_{i=2}^{k}c_i \lambda_i \\
\lambda_2 & \geq &  \frac{n_2}{d_2}\lambda_3 \nonumber\\ \nonumber
\lambda_3 & \geq &  \frac{n_3}{d_3}\lambda_4 \\ \nonumber
& \vdots & \\ \nonumber
\lambda_{k-1} & \geq &  \frac{n_{k-1}}{d_{k-1}}\lambda_k \\ \nonumber
\lambda_{k} & \geq & 0   \\ \nonumber
\label{rational}
\end{eqnarray}
where:
\begin{itemize}
\item For $1 \leq i \leq k-1$, $n_i$ and $d_i$ are positive integers and
\item
the $c_i$ are any integers which make the first constraint 
strong enough to guarantee that $\lambda_1 \geq 0$.
(We will see several examples of such $c_i$.)
\end{itemize}
For $i = 1 \ldots k$, let
\[
a_i = \prod_{j=1}^{i-1}d_j \prod_{t=i}^{k-1}n_t.
\]
Then $a_i/a_{i+1} = n_i/d_i$, so the system (\ref{rational})
 above is equivalent to:
\begin{eqnarray}
\lambda_1 & \geq & c_1 \lceil \frac{a_1}{a_2}\lambda_2 \rceil +
\sum_{i=2}^{k}c_i \lambda_i \ \ \ \mbox{and}
\label{lhp2}
\end{eqnarray}
\[
\lambda_2/a_2 \geq
\lambda_3/a_3 \geq \cdots \geq
\lambda_{k-1}/a_{k-1} \geq
\lambda_k/a_k.
\]

\subsection{The Generating Function for (\ref{rational})}
Consider first the case where $c_1=1$ and $c_i=0$ for $i \geq 2$.

\begin{theorem}
Given a sequence of positive integers
$a_1, \ldots, a_k$,
the generating function for the compositions
$\lambda_1, \ldots, \lambda_k$ satisfying
\begin{eqnarray*}
\lambda_1 &  \geq & a_1 \lambda_2/a_2\\
\lambda_2 &  \geq & a_2 \lambda_3/a_3\\
& \vdots & \\
\lambda_{k-2} &  \geq & a_{k-2} \lambda_{k-1}/a_{k-1}\\
\lambda_{k-1} &  \geq & a_{k-1} \lambda_{k}/a_{k}\\
\lambda_{k} &  \geq & 0
\end{eqnarray*}
is
\begin{equation}
\frac{
 \sum_{z_2=0}^{a_2-1}
\sum_{z_3=0}^{a_3-1} \cdots
\sum_{z_{k}=0}^{a_{k}-1} 
q^{\lceil \frac{a_1z_2}{a_2} \rceil + \sum_{i=2}^{k}z_i}
\prod_{i=2}^{k-1}q^{b_i\lceil \frac{z_{i+1}}{a_{i+1}} - \frac{z_{i}}{a_{i}}
\rceil}}
{\prod_{i=1}^{k} (1-q^{b_i})}
\label{gfnoc}
\end{equation}
where $b_1=1$ and for $2 \leq i \leq k$,
\[ b_i = a_1 + a_2 +\cdots + a_{i}.
\]
\end{theorem}

\noindent
{\bf Proof}
For $2 \leq i \leq k$, let
\[
\lambda_i = a_ix_i + z_i,
\]
where $x_i \geq 0 $ and $0 \leq z_i < a_i$.
The system of inequalities becomes
\begin{eqnarray*}
\lambda_1 &  \geq & a_1 x_2 + a_1 z_2/a_2\\
a_2x_2 + z_2 &  \geq & a_2 x_3 + a_2 z_3/a_3\\
a_3x_3 + z_3 &  \geq & a_3 x_4 + a_3 z_4/a_4\\
& \vdots & \\
a_{k-2}x_{k-2}+z_{k-2} & \geq & a_{k-2} x_{k-1}+a_{k-2}z_{k-1}/a_{k-1}\\
a_{k-1}x_{k-1}+z_{k-1} & \geq & a_{k-1} \lambda_{k}\\
a_{k}x_{k}+z_{k} & \geq & 0 \\
\end{eqnarray*}
Rearrange the  sequence              
$\lambda_1,  \lambda_2, \ldots, \lambda_k $
by moving $(a_i-1)x_i + z_i$ ``dots" in the Ferrers diagram
 from part $\lambda_i$ to part $\lambda_1$ 
for
$2 \leq i \leq k$ to get a new sequence of the same weight:
\[
x_1, x_2, \ldots , x_{k},
\]
where
\[
x_1 = \lambda_1 +  \sum_{i=2}^{k}((a_i-1)x_i+z_i)
\]
satisfying: 
\begin{eqnarray}
x_1 &  \geq &  (a_1+a_2-1)x_2 + \sum_{i=3}^{k}(a_i-1)x_i + a_1z_2/a_2 +
\sum_{i=2}^{k}z_i  \nonumber\\
x_2  &  \geq & x_3 +  z_3/a_3 - z_2/a_2 \nonumber \\
x_3  &  \geq & x_4 +  z_4/a_4 - z_3/a_3 \nonumber \\
& \vdots &  \nonumber \\
x_{k-2} & \geq & x_{k-1}+z_{k-1}/a_{k-1}-z_{k-2}/a_{k-2} \nonumber \\
x_{k-1} & \geq & x_{k}+z_{k}/a_{k}-z_{k-1}/a_{k-1} \nonumber \\
x_{k} & \geq & 0.   \nonumber \\
\label{lhp}
\end{eqnarray}    
This is just the system
\begin{eqnarray*}
x_1 &  \geq &  (a_1+a_2-1)x_2 + \sum_{i=3}^{k}(a_i-1)x_i + s_1\\
x_i  &  \geq & x_{i+1} + s_{i}\ \ \ 2\le i\le k-1\\
x_k & \geq & 0
\end{eqnarray*}
where
\[
s_1 = \lceil a_1z_2/a_2 \rceil + \sum_{i=2}^{k}z_i
\]
and
\[
s_i = \lceil z_{i+1}/a_{i+1} - z_{i}/a_{i} \rceil
\]
for $2 \leq i \leq k$.  

Thus, by Corollary 3 and Corollary 4,
the generating function for
fixed $s_1, \dots s_{k-1}$ is
\[
\frac{
\prod_{i=1}^{k-1}q^{b_is_i}}
{\prod_{i=1}^{k} (1-q^{b_i})},
\]
where $b_1=1$ and for $2 \leq i \leq k$,
$ b_i = a_1 + a_2 +\cdots + a_{i}.$
Summing over all possible sequences
$s_1, s_2, \ldots, s_{k-1}$ as the $z_i$ vary independently from
$0$ to $a_i-1$ gives the result.
\hfill $\Box$

\noindent
{\bf Example 14} Suppose $(a_1,a_2,a_3,a_4)=(4,3,2,1)$.
Then $(b_1,b_2,b_3,b_4)=(1,7,9,10)$
so the denominator of (\ref{gfnoc}) is the product
$(1-q)(1-q^7)(1-q^9)(1-q^{10})$.
The numerator of (\ref{gfnoc}) is the sum,
as $(z_2,z_3,z_4)$ range over the set
$\{
(0,0,0),
(1,0,0),
(2,0,0),
(0,1,0),$
$(1,1,0),
(2,1,0) \}$,
 of the terms
\[
q^{\lceil 4z_2/3 \rceil + z_2+z_3+z_4} q^{7 \lceil z_3/2-z_2/3 \rceil}
q^{9 \lceil z_4/1-z_3/2 \rceil} =
q^{\lceil 4z_2/3 \rceil + z_2+z_3} q^{7 \lceil z_3/2-z_2/3 \rceil}.
\]
So, the generating function is
\begin{eqnarray*}
\frac{1+q^3+q^5+q^8+q^{11}+q^6}
{(1-q)(1-q^7)(1-q^9)(1-q^{10})}
&=&
\frac{(1+q^3+q^6)(1+q^5)}
{(1-q)(1-q^7)(1-q^9)(1-q^{10})}\\
&=&
\frac{1}
{(1-q)(1-q^7)(1-q^3)(1-q^{5})}.
\end{eqnarray*}
More generally,
in the special case of Theorem 2 that $a_i = k-i+1$ for $1 \leq i \leq k$ and
where $c_1=1$  and $c_i=0$ for $1=2, \ldots k$ ,
the constraints become:
\[
\lambda_1/k \geq
\lambda_2/(k-1) \geq \cdots \geq
\lambda_{k-1}/2 \geq
\lambda_k/1
\]
and the sequences satisfying these constraints are partitions
known as {\em lecture hall partitions}.
It was first shown in \cite{BME1} that these partitions have
generating function
\[\prod_{i=1}^{k}\frac{1}{(1-q^{2i-1})}, \]
giving a finite form of the identity of Euler equating partitions of $n$
into distinct parts with partitions of $n$ into odd parts.

Note that although Theorem 2 gives an explicit form of the generating
function
for these partitions, it does not help with the factoring of the
numerator.   In another proof of the lecture hall partitions
theorem 
\cite{Om1}, Andrews uses partition analysis to get the generating
function and shows how to factor the numerator
via a permutation of the set of tuples $(z_2, \ldots z_k)$.

\noindent
{\bf Example 15}
Suppose $(a_1,a_2,a_3,a_4,a_5)=(1,3,2,3,1)$.
Then $(b_1,b_2,b_3,b_4,b_5)=(1,4,6,9,10)$
so the denominator of (\ref{gfnoc}) is the product
$(1-q)(1-q^4)(1-q^6)(1-q^9)(1-q^{10})$.
The numerator of (\ref{gfnoc}) is the sum,
as $(z_2,z_3,z_4,z_5)$ ranges over the set
$\{
(0,0,0,0),$
$(1,0,0,0),$
$(2,0,0,0),$
$(0,1,0,0),$
$(1,1,0,0),$
$(2,1,0,0),$
$(0,0,1,0),$
$(1,0,1,0),$
$(2,0,1,0),$
$(0,1,1,0),$
$(1,1,1,0),$
$(2,1,1,0),$
$(0,0,2,0),$
$(1,0,2,0),$
$(2,0,2,0),$
$(0,1,2,0),$
$(1,1,2,0),$
$(2,1,2,0) \}$,
 of the terms
{
\begin{eqnarray*}
 q^{\lceil z_2/3 \rceil + z_2+z_3+z_4} q^{4 \lceil z_3/2-z_2/3 \rceil}
q^{6 \lceil z_4/3-z_3/2 \rceil} 
q^{9 \lceil z_5/1-z_4/3 \rceil} .
\end{eqnarray*}}
So, the generating function is
{\small\begin{eqnarray*}
 \frac{1+q^2+q^3+q^5+q^{7}+q^4+
q^7 + q^9 + q^{10} + q^6 + q^8 + q^5 +
q^8 + q^{10} + q^{11} + q^{13} + q^{15} + q^{12}}
{(1-q)(1-q^4)(1-q^6)(1-q^9)(1-q^{10})}\\
=
\frac{(1+q^2+q^4)(1+q^3+q^6)(1+q^5)}
{(1-q)(1-q^4)(1-q^6)(1-q^9)(1-q^{10})}\\
=
\frac{1}
{(1-q)(1-q^2)(1-q^3)(1-q^4)(1-q^{5})}.
\end{eqnarray*}}
Thus there is a one-to-one correspondence between the 
compositions of $n$ into 5 parts satisfying
\[
\frac{\lambda_1}{1} \geq
\frac{\lambda_2}{3} \geq
\frac{\lambda_3}{2} \geq
\frac{\lambda_4}{3} \geq
\frac{\lambda_5}{1}
\]
and the partitions of $n$ into parts of size at most 5.

It is straightforward to extend Theorem 2 to get the generating
function for the system (\ref{lhp2}) with more general $c_i$.

\begin{corollary}
Given   a sequence of positive integers
$a_1, \ldots, a_k$ and a sequence of integers
$c_1, \ldots, c_k$, consider the set of
sequences
$\lambda_1, \ldots, \lambda_k$ satisfying
\begin{eqnarray}
\lambda_1 & \geq & c_1 \lceil \frac{a_1}{a_2}\lambda_2 \rceil +
\sum_{i=2}^{k}c_i \lambda_i \ \ \ {\mbox and}
\nonumber
\end{eqnarray}
\[
\lambda_2/a_2 \geq
\lambda_3/a_3 \geq \cdots \geq
\lambda_{k-1}/a_{k-1} \geq
\lambda_k/a_k.
\]
As long as
the $c_i$ are integers which make the first constraint 
strong enough to guarantee that $\lambda_1 \geq 0$, the generating
function 
is
\begin{equation}
\frac{
 \sum_{z_2=0}^{a_2-1}
\sum_{z_3=0}^{a_3-1} \cdots
\sum_{z_{k}=0}^{a_{k}-1} 
q^{c_1\lceil \frac{a_1z_2}{a_2} \rceil + \sum_{i=2}^{k}(c_i+1)z_i}
\prod_{i=2}^{k-1}q^{b_i\lceil \frac{z_{i+1}}{a_{i+1}} - \frac{z_{i}}{a_{i}}
\rceil}}
{\prod_{i=1}^{k} (1-q^{b_i})},
\label{gfc}
\end{equation}
where $b_1=1$, $b_2=c_1a_1$ and for $3 \leq i \leq k$,
\[ b_i = c_1a_1 + (c_2+1)a_2 +\cdots + (c_{i}+1)a_{i}.
\]
\end{corollary}

\noindent
{\bf Proof}
Using the same strategy as in the proof of Theorem 2, the system of
inequalities in
the Corollary becomes
\begin{eqnarray*}
x_1 &  \geq &  (c_1a_1+(c_2+1)a_2-1)x_2 + \sum_{i=3}^{k}((c_i+1)a_i-1)x_i + s_1\\
x_i  &  \geq & x_{i+1} + s_{i}\ \  \ 2\le i\le k-1\\
x_k & \geq & 0
\end{eqnarray*}
where now $0 \leq z_i \leq i-1$,
\[
s_1 = c_1\lceil a_1z_2/a_2 \rceil + \sum_{i=2}^{k}(c_i+1)z_i,
\]
and
\[
s_i = \lceil z_{i+1}/a_{i+1} - z_{i}/a_{i} \rceil
\]
for $2 \leq i \leq k-1$.  
Thus, by Corollary 3 and Corollary 4,
the generating function for
fixed $s_1, s_2, \ldots s_{k-1}$ is
\[
\frac{
\prod_{i=1}^{k-1}q^{b_is_i}}
{\prod_{i=1}^{k} (1-q^{b_i})},
\]
where $b_1=1$, $b_2=c_1a_1$ and for $3 \leq i \leq k$,
$ b_i = c_1a_1 + (c_2+1)a_2  + (c_3+1)a_3+\cdots +(c_{i}+1) a_{i}.$
Summing over all possible sequences
$s_1, s_2, \ldots s_{k-1}$ as the $z_i$ vary independently from
$0$ to $a_i-1$ gives the result.
\hfill $\Box$


\noindent
{\bf Example 16}
Consider sequences $(\lambda_1,\lambda_2,\lambda_3,\lambda_4,\lambda_5)$
satisfying 
\[\lambda_1 \geq  \lambda_2/2  +3 \lambda_5
\]
and
\[\frac{\lambda_2}{2} \geq \frac{\lambda_3}{1} \geq 
\frac{\lambda_4}{2} \geq \frac{\lambda_5}{1} \geq 0.
\]
The constraints guarantee these $\lambda$ are compositions and they
satisfy the conditions of Corollary 5 with
$k=5$,
$(a_1,a_2,a_3,a_4,a_5)=(1,2,1,2,1)$ and
$(c_1,c_2,c_3,c_4,c_5)=(1,0,0,0,3)$.
Then $(b_1,b_2,b_3,b_4,b_5)=(1,3,4,6,10)$
so the denominator of (\ref{gfc}) is the product
$(1-q)(1-q^3)(1-q^4)(1-q^6)(1-q^{10})$.
The numerator of (\ref{gfc}) 
 is the sum of the terms
\[
q^{ \lceil z_2/2 \rceil  + z_2+z_3+z_4+4z_4} 
q^{3 \lceil z_3/1-z_2/2 \rceil}
q^{4 \lceil z_4/2-z_3/1 \rceil}
\]
as $(z_2,z_3,z_4,z_5)$ ranges over the set
$\{
(0,0,0,0),
(1,0,0,0),
(0,0,1,0),
(1,0,1,0),
 \}$.
So, by Corollary 5, the generating function is
\begin{eqnarray*}
\frac{1+q^2+q^5+q^7}
{(1-q)(1-q^3)(1-q^4)(1-q^6)(1-q^{10})}
&=&
\frac{(1+q^2)(1+q^5)}
{(1-q)(1-q^3)(1-q^4)(1-q^6)(1-q^{10})}\\
&=&
\frac{1}
{(1-q)(1-q^2)(1-q^3)(1-q^{5})(1-q^6)}.
\end{eqnarray*}

\noindent
{\bf Example 17}
Consider sequences $(\lambda_1,\lambda_2,\lambda_3,\lambda_4)$
satisfying 
\[\lambda_1 \geq 2 \lceil  7\lambda_2/3 \rceil +  3 \lambda_2+ \lambda_3 +5 \lambda_4
\]
and
\[\frac{\lambda_2}{3} \geq \frac{\lambda_3}{2} \geq 
\frac{\lambda_4}{1} \geq 0.
\]
The constraints guarantee these $\lambda$ are partitions and they
satisfy the conditions of Corollary 5 with
$k=4$,
$(a_1,a_2,a_3,a_4)=(7,3,2,1)$ and
$(c_1,c_2,c_3,c_4)=(2,3,1,5)$.
Then $(b_1,b_2,b_3,b_4)=(1,26,30,36)$
so the denominator of (\ref{gfc}) is the product
$(1-q)(1-q^{26})(1-q^{30})(1-q^{36})$.
The numerator of (\ref{gfc}) 
 is the sum of the terms
\[
q^{ 2\lceil 7z_2/3 \rceil  + 4z_2+2z_3 +6z_4} q^{26 \lceil z_3/2-z_2/3 \rceil}
 q^{30 \lceil z_4/1-z_3/2 \rceil}
\]
as $(z_2,z_3,z_4)$ ranges over the set
$\{
(0,0,0),
(1,0,0),
(2,0,0),
(0,1,0),
(1,1,0),
(2,1,0) \}$.
So, by Corollary 5, the generating function is
\begin{eqnarray*}
\frac{1+q^{10}+q^{18}+q^{28}+q^{38}+q^{20}}
{(1-q)(1-q^{26})(1-q^{30})(1-q^{36})}
&=&
\frac{(1+q^{10}+q^{20})(1+q^{18})}
{(1-q)(1-q^{26})(1-q^{30})(1-q^{36})}\\
&=&
\frac{1}
{(1-q)(1-q^{10})(1-q^{18})(1-q^{26})}.
\end{eqnarray*}

\subsection{A Special Case}

\begin{corollary}
Suppose the sequence $a_1, a_2, \ldots, a_k$ has the property that
for $1 \leq i \leq k-1$ if $a_i > 1$, then $a_{i+1}=1$.
Then there is a one-to-one correspondence between the compositions
$\lambda_1, \ldots, \lambda_k$ of $n$ satisfying
\[
\lambda_1/a_1 \geq
\lambda_2/a_2 \geq
\lambda_3/a_3 \geq \cdots \geq
\lambda_{k-1}/a_{k-1} \geq
\lambda_k/a_k
\]
and
the partitions of $n$ into parts
in
\[
\{1, b_2, b_2+1, \ldots, b_{k}\},
\]
where $b_1=1$,
$b_2 = a_1+a_2$ and $b_{i+1} = b_{i}+a_{i+1} $ for $2 \leq i \leq k-1$,
such that at most one part can appear from each of the sets
\[
S_i = \{p | b_i+1 \leq p \leq b_{i+1}-1\}.
\]
\end{corollary}

\noindent
{\bf Proof}
Since $z_j = 0$ whenever $a_j=1$, the generating function
(\ref{gfnoc}) becomes

\begin{equation*}
\frac{
 \sum_{z_2=0}^{a_2-1}
\sum_{z_3=0}^{a_3-1} \cdots
\sum_{z_{k}=0}^{a_{k}-1} 
q^{\lceil \frac{a_1z_2}{a_2} \rceil + z_2}
\prod_{i=2: a_{i+1} > 1 }^{k-1}q^{b_i\lceil \frac{z_{i+1}}{a_{i+1}}
\rceil + z_{i+1}}}
{\prod_{i=1}^{k} (1-q^{b_i})}
\end{equation*}
and since consecutive $a_i$ cannot both be greater than 1, we get
\begin{equation*}
\frac{
 \sum_{z_2=0}^{a_2-1}
q^{\lceil \frac{a_1z_2}{a_2} \rceil + z_2}
\prod_{i=2: a_{i+1} > 1 }^{k-1}
\sum_{z_{i+1}=0}^{a_{i+1}-1} 
q^{b_i\lceil \frac{z_{i+1}}{a_{i+1}}
\rceil + z_{i+1}}}
{\prod_{i=1}^{k} (1-q^{b_i})}
\end{equation*}  
which, letting $b_1=1$, gives
\begin{equation*}
\frac{
\prod_{i=1: a_{i+1} > 1 }^{k-1}
\sum_{z_{i+1}=0}^{a_{i+1}-1} 
q^{b_i\lceil \frac{z_{i+1}}{a_{i+1}}
\rceil + z_{i+1}}}
{\prod_{i=1}^{k} (1-q^{b_i})}.
\end{equation*}   
Since $b_2 = a_1+a_2$ and $b_{i+1} = b_{i}+a_{i+1}$
we get
\begin{equation}
\frac{
\prod_{i=1: b_{i+1}-b_i > 1 }^{k-2}
(1 + q^{b_i+1} + q^{b_i+2} + \cdots + q^{b_{i+1}-1})}
{\prod_{i=1}^{k} (1-q^{b_i})}.
\label{1212gf}
\end{equation}   
Note that each integer $m$ satisfying
$b_2 \leq m \leq b_{k}$ appears exactly once as an exponent
in either the numerator or denominator.
So (\ref{1212gf}) is the generating function for partitions into parts
in
\[
\{1, b_2, b_2+1, \ldots, b_{k}\}
\]
such that at most one part can appear from each of the sets
\[
S_i = \{p\ |\ b_i+1 \leq p \leq b_{i+1}-1\}.
\]
(The set $S_i$ is empty
if
$b_{i+1}-b_i > 1 $). This result can easily be proved bijectively. 
\hfill $\Box$













\noindent
{\bf Example 18}
Compositions of $n$ satisfying
\[
\lambda_1 \geq \lambda_2 \geq \frac{\lambda_3}{2} \geq
 \lambda_4 \geq \frac{\lambda_5}{2}
\geq \lambda_6 \geq  \cdots \geq \lambda_{2k} \geq \frac{\lambda_{2k+1}}{2} \geq 0
\]
are in one-to-one correspondence with the set of partitions of $n$
into parts of size at most $3k+2$ in which parts divisible by 3 can
appear at most once.

\noindent
{\bf Example 19} Compositions of $n$ satisfying
\[
\lambda_1 \geq \frac{\lambda_2}{2} \geq \lambda_3 \geq 
\frac{\lambda_4}{3} \geq \lambda_5
\geq \frac{\lambda_6}{4} \geq
 \cdots \geq \lambda_{2k-1} \geq \frac{\lambda_{2k}}{k+1}
\]
are in one-to-one correspondence with the set of partitions of $n$
into parts of size at most $k(k+5)/2$, such that for each $j \geq 3$, at 
most one part can occur from the set
\[
\left\{i \ | \ \binom{j}{2} \leq i+1 \leq \binom{j}{2} + (j-3) \right\}.
\]

\section{Two Variable Generating Functions}


In their study of lecture hall partitions,  Bousquet-M\'elou and
Eriksson found it very useful to consider the 2-variable
(odd/even weighted) generating function of the set of partitions
satisfying the lecture hall constraints.
We show here how our method can be adapted to get 
multivariable generating functions for compositions satisfying
linear constraints, using the two-variable case as an example.

Given a sequence $\lambda=(\lambda_1,\dots, \lambda_k)$, we denote
by $\lambda_o$  the subsequence 
 $(\lambda_1,\lambda_3,\lambda_5,\dots)$ 
and by $\lambda_e$  the subsequence 
 $(\lambda_2,\lambda_4,\lambda_6,\dots)$.


\subsection{Integer Coefficients}

Let $A[1..k,1..k]$ be a strictly upper triangular matrix of 
integers,   such that $P_A$ is a set of compositions.
Let $Q_A(l,m,k)$ be the subset of $P_A$ consisting of those 
compositions 
$\lambda=(\lambda_1,\lambda_2, \cdots, \lambda_k)$, 
with $l=|\lambda_o|= \lambda_1 + \lambda_3 + \lambda_5 + \cdots$ and
with $m=|\lambda_e|= \lambda_2 + \lambda_4 + \lambda_6 + \cdots$,
satisfying
\[
\lambda_i \geq \sum_{j=i+1}^{k} A[i,j] \lambda_{j} \ \ \ 
{\mbox{\rm  for $i \geq 1$}}.
\]
We would now like to write the generating function
\[
\sum_{l,m \geq 0}|Q_A(l,m,k)|x^ly^m.
\]

\begin{theorem}
Let $A$ a $k\times k$ strictly upper triangular matrix of integers,. 
If $B=I+A+A^2+\ldots +A^{k-1}$ has nonnegative coefficients
then the two-variable generating function for $Q_A(l,m,k)$ is
\begin{eqnarray}
\prod_{i=1}^{k}\frac{1}{1-x^{o_i}y^{e_i}},
\label{gf2}
\end{eqnarray}
where $o_i =\sum_{j \geq 1} B[2j-1,i]$ and $e_i =\sum_{j \geq 1} B[2j,i]$.
\end{theorem}
\noindent
{\bf Proof}
As $\lambda=B[s_1,\ldots, s_k]^T$ then
$\lambda_{2i+1}=\sum_{j}B[2j+1,j]s_j$ and
$\lambda_{2i}=\sum_{j}B[2j,i]s_j$. Therefore
$|\lambda_o|=\sum_{i\ge 1}o_is_i$ and
$|\lambda_e|=\sum_{i\ge 1}e_is_i$. The result follows.
\hfill $\Box$

\noindent
{\bf Example 20}
Let $f_k(x,y)$ be the generating function for ordinary partitions:
\[
\lambda_1 \geq \lambda_2 \geq  \cdots \geq \lambda_k \geq 0.
\]
Then 
$o_1, o_2, o_3, \ldots$ is: $1,1,2,2,3,3,4,4,\ldots$ 
and $e_1, e_2, e_3, \ldots$ is: $0,1,1,2,2,3,3,4,4,\ldots$ 
so the generating function is
\[
f_k(x,y) =
\prod_{i=1}^{\lceil k/2 \rceil}\frac{1}{1-x^iy^{i-1}}
\prod_{i=1}^{\lfloor k/2 \rfloor}\frac{1}{1-x^iy^{i}}.
\]
Note that 
\[
f_{k+1}(x,y) = \frac{f_k(x^2y,x^{-1})}{(1-x)}.
\]


\subsection{Alpha-beta Sequences}

Let $A$ be a $k \times k$ strictly upper triangular matrix of
integers, 
 such that $P_A$ is a set of {\em partitions} and let $G(x,y)$
 be the two variable generating function, that is,
\begin{eqnarray}
G(x,y) = \sum_{l,m \geq 0}|Q_A(l,m)|x^ly^m.
\label{Gxy}
\end{eqnarray}

\begin{theorem}
Let $\alpha \geq 1$ and $\beta \leq \alpha$ be integers (note $\beta$ can be
negative).
The set of compositions $\lambda_0, \lambda_1, \ldots \lambda_k$,
satisfying
\[
\lambda_0 = \alpha(\lambda_{1} + \lambda_{3} + \lambda_{5} + \cdots) 
	- \beta(\lambda_{2} + \lambda_{4} + \lambda_{6} + \cdots)
\]
and for $1 \leq i \leq k$, 
$\lambda_i \geq \sum_{j=i+1}^{k}A[i,j]\lambda_{j}$
has generating function 
$G(x^{\alpha}y,x^{-\beta + 1})$, where
$G$ is defined by (\ref{Gxy}).
\end{theorem}

\noindent
{\bf Proof}  Let partition $\lambda =
(\lambda_1, \lambda_2, \ldots \lambda_k)$,
 satisfy the  constraints
$\lambda_i \geq \sum_{j=i+1}^{k}A[i,j]\lambda_{j}$.
  Prepend part
$\lambda_0$ to satisfy the new constraint and call the new composition
$\lambda'$.
(Since $\lambda$ is a partition, the conditions on $\alpha,\beta$
ensure $\lambda_0 \geq 0$.
The  odd parts of $\lambda$ become the even parts of
$\lambda'$.
The odd parts of $\lambda'$ are the
 even parts of $\lambda$ plus the new part $\lambda_0$, so
\[|\lambda'|_e = |\lambda|_o;  \ \ \
|\lambda'|_o = |\lambda|_e + \alpha |\lambda|_o -\beta |\lambda|_e.\ \ \ \ \ \ \ \ \ \ \ \ \ \ \Box
\]

\begin{corollary}
For fixed constants $\alpha \geq 1$ and $\beta \leq \alpha$,
the sequences
$(\lambda_1, \ldots \lambda_k)$
satisfying, for $1 \leq i \leq k$,
\[
\lambda_i \geq \alpha(\lambda_{i+1} + \lambda_{i+3} + \lambda_{i+5} +
\cdots) 
	- \beta(\lambda_{i+2} + \lambda_{i+4} + \lambda_{i+6} + \cdots),
\]
form a set of compositions with odd/even generating function
\[
G_k(x,y) = \prod_{i=1}^{k}\frac{1}{(1-x^{o_i}y^{e_i})}
\]
where
\[
o_1=1;\ \ o_2=\alpha;  \ \    o_i = \alpha o_{i-1} + (-\beta + 1)o_{i-2}
\ \ \mbox{for $i>2$}
\]
and
\[
e_1=0; \ \  e_2=1;  \ \  e_i = o_{i-1} 
\ \ \mbox{for $i>2$}.
\]
\label{alphabeta}
\end{corollary}

\noindent
{\bf Proof}
The constraints ensure  that $\lambda_i \geq 0$ for $i=1, \ldots k$, so as
in proof of
Theorem 4
their generating functions satisfy the recurrence
\[ G_{k+1}(x,y) = G_k(x^{\alpha}y,x^{-\beta + 1})/(1-x).\ \ \ \ \ \ \ \ \ \ \ \ \ \Box
\]


\noindent {\bf Example 21}
If $\alpha=1$ and $\beta = -1$  in Corollary \ref{alphabeta}
we get:
\[
o_1=1;\ \ o_2=1; \ \   o_i = o_{i-1}+2o_{i-2}\ \
\]
and
\[
e_0=1; \ \  e_1=2; \ \ e_i = o_{i-1},
\] so the generating function is
\[
\frac{1}{(1-x)(1-xy)(1-x^3y)(1-x^5y^3)(1-x^{11}y^5)(1-x^{21}y^{11})
\cdots}.
\]
Substituting $x=y=q$ tells us that the partitions 
$\lambda_1, \ldots,  \lambda_k$ of $n$ satisfying
\[ \lambda_i \geq \lambda_{i+1}+\lambda_{i+2}+ \lambda_{i+3} + \cdots
\]
are in one-to-one correspondence with partitions of $n$
whose parts must be one of the first $n$ powers of 2: $1, 2, \ldots, 2^{n-1}$.

\noindent {\bf Example 22}
In Corollary \ref{alphabeta},
if $\alpha=1$ and $\beta = 1$, we get:
\[
o_1=1;\ \ o_2=1; \ \   o_i = o_{i-1}\ \
\]
and
\[
e_1=0; \ \  e_2=1; \ \ e_i = o_{i-1},
\] so the generating function is
\[
\frac{1}{(1-x)(1-xy)(1-xy)(1-xy)(1-xy)(1-xy)
\cdots}.
\]
Substituting $x=y=q$ tells us that the compositions
$\lambda_1, \ldots,  \lambda_k$ of $n$ satisfying
\[ \lambda_i \ge \sum_{j=1}^{k-i} (-1)^{j} \lambda_{i+j}
\]
are in one-to-one correspondence with 
{\em compositions} of $n$
into $k$ parts whose first part can be any nonnegative integer but
whose other parts must be nonnegative even integers.
(Compare with Example 2.)

\noindent {\bf Example 23}
If $\alpha=2$ and $\beta = 2$, from Corollary \ref{alphabeta} we get:
\[
o_1=1;\ \ o_2=2; \ \   o_i = 2o_{i-1}-o_{i-2}\ \
\]
and
\[
e_1=0; \ \  e_2=1; \ \ e_i = o_{i-1},
\] so the generating function is
\[\frac{1}{(1-x)(1-x^2y)(1-x^3y^2)(1-x^4y^3)(1-x^{5}y^4)(1-x^{6}y^{5})
\cdots}.
\]
Substituting $x=y=q$ tells us that the partitions 
$\lambda_1, \ldots,  \lambda_k$ of $n$ satisfying
\[ \lambda_i \ge 2 \sum_{j=1}^{k-i} (-1)^{j} \lambda_{i+j}
\]
are in one-to-one correspondence with 
partitions into odd parts less than $2k$.


\subsection{The Two Variable Generating Function for Rational
Diagonal Constraints}

To get the odd/even weighted generating function for the
compositions satisfying the constraints:
\begin{equation}
\frac{\lambda_1}{a_{1}} \geq
\frac{\lambda_2}{a_{2}} \geq \frac{\lambda_3}{a_{3}} \geq \cdots \geq
\frac{\lambda_{k-1}}{a_{k-1}} \geq \frac{\lambda_k}{a_k} \geq 0,
\label{lhpconstraints}
\end{equation}
in the system (\ref{lhp}) we separate the ``dots" comprising $x_1$
into
those which come from odd parts:
\[
a_1x_2+(a_3-1)x_3+(a_5-1)x_5+ \cdots + z_3+z_5+ \cdots
\]
and those which come from even parts:
\[
(a_2-1)x_2+(a_4-1)x_4+ \cdots + z_2+z_4+z_6 + \cdots
\]
Then we get the following refinement of
(\ref{gfnoc}):

\begin{equation}
G(x,y) =
\frac{
 \sum_{z_2=0}^{a_2-1}
\cdots
\sum_{z_{k}=0}^{a_{k}-1} 
x^{\lceil \frac{a_1z_2}{a_2} \rceil + z_3 + z_5 + \cdots }
y^{z_2+z_4   + \cdots}
\prod_{i=2}^{k-1}(x^{o_i}y^{e_i})^{\lceil \frac{z_{i+1}}{a_{i+1}} - \frac{z_{i}}{a_{i}} \rceil}}
{\prod_{i=1}^{k} (1-x^{o_i}y^{e_i})},
\label{2gf}
\end{equation}
where
$o_1 = 1$, $e_1=0$, and for $2 \leq i \leq k$
$o_i = a_1 + a_3 + a_5 + \cdots + a_{2 \lfloor (i-2)/2 \rfloor + 1}$ and
$e_i = a_2 + a_4 + a_6 + \cdots + a_{2 \lceil (i-1)/2 \rceil }$.

\subsection{Variations on Lecture Hall Partitions}

\begin{theorem}
If $G(x,y)=H(x,y)/(1-x)$ is the generating function given in (\ref{2gf}), then
whenever $a_1 \geq a_2 \geq \cdots \geq a_k$,
for any $l > 0$ and $j\ge 2-l$,
$H(x^l,x^{j-1}y)$ is the generating function for the partitions
satisfying
\begin{eqnarray}
\lambda_1  &  = & l \lceil a_1 \lambda_2/a_2 \rceil +
\sum_{i\ge 1}((j-1)\lambda_{2i} 
+(l-1)\lambda_{2i+1})  \nonumber\\
\frac{\lambda_2}{a_{2}} &  \geq  & \frac{\lambda_3}{a_{3}} \geq \cdots \geq
\frac{\lambda_{k-1}}{a_{k-1}} \geq \frac{\lambda_k}{a_k} \geq 0.
\label{lhpgenconstraints}
\end{eqnarray}
\label{lhpgen}
\end{theorem}      

\noindent
{\bf Proof}
If $j=l=1$, the system (\ref{lhpgenconstraints}) is the same as
(\ref{lhpconstraints}), except with equality for $\lambda_1$,
 and the generating function is $H(x,y)=G(x,y)/(1-x)$.
Suppose $\lambda$ satisfies (\ref{lhpconstraints}).
To transform $\lambda$ into a composition satisfying
(\ref{lhpgenconstraints}), we increase the first part by
$(l-1)|\lambda|_o +
(j-1)|\lambda|_e$ to get $\lambda'$.
The conditions of the theorem guarantee that this increase is
nonnegative.
Then 
$|\lambda'|_e = |\lambda|_e$ and
\[
|\lambda'|_o =
|\lambda|_o  +
(l-1)|\lambda|_o + (j-1)|\lambda|_e =
l|\lambda|_o + (j-1)|\lambda|_e.\ \ \ \ \ \ \ \ \ \ \ \ \ \ \ \ \ \ \ \ \Box
\]

\noindent{\bf Remark} If the $a_i$s are not non-increasing then the theorem
is still true if $l$ and $j$ are both positive.



We can use Theorem \ref{lhpgen} to generalize the Lecture Hall Partition
theorem of Bousquet-M\'elou and Eriksson.
\begin{corollary}
For $l>0$ and $j\ge 2-l$,
the generating function for the sequences
$\lambda_1, \ldots, \lambda_k$ satisfying
\[
\lambda_1   \geq  l \lceil k \lambda_2/(k-1) \rceil +
(j-1)(\lambda_2 + \lambda_4 + \lambda_6 + \cdots) +
(l-1)(\lambda_3 + \lambda_5 + \lambda_7 + \cdots) 
\]
and
\[
\frac{\lambda_2}{k-1} \geq \frac{\lambda_3}{k-2} \geq
\cdots \geq
\frac{\lambda_{k-1}}{2} \geq \frac{\lambda_k}{1}\geq 0 
\]
is
\[
\frac{1}{(1-q)}\prod_{i=1}^{k-1}\frac{1}{(1-q^{il+ij+l})}.
\]
\end{corollary}

\noindent
{\bf Proof}
 Bousquet-M\'elou and Eriksson have shown  in \cite{BME1} that
when $j=l=1$ the generating function is:
\[
G(x,y) = \prod_{i=0}^{k-1}\frac{1}{(1-x^{i+1}y^i)}.
\]
The result follows from Theorem \ref{lhpgen},
setting $x=y=q$. \hfill $\Box$


\noindent
{\bf Example 24}
\[
\prod_{i=1}^{k-1}\frac{1}{(1-q^{2i+3})}
\]
is the generating function for the partitions satisfying
\[
\lambda_1   =  3 \lceil k \lambda_2/(k-1) \rceil +
2 \sum_{j=1}^{k-1}(-1)^j \lambda_{1+j}
\]
and
\[
\frac{\lambda_2}{k-1} \geq \frac{\lambda_3}{k-2} \geq
\cdots \geq
\frac{\lambda_{k-1}}{2} \geq \frac{\lambda_k}{1} .
\]
This is Corollary 8 with $l=3$ and $j=-1$.

\section{Conclusion}

In this paper we have shown that  compositions
and partitions defined by (in)equalities have nice and easy-to-compute
generating functions in the integer case and, sometimes,
in the rational case. The proofs are all bijective and the techniques
are powerful and well suited.
However these techniques 
are not completely satisfactory for analyzing families like the
Lecture Hall partitions, since additional (and non-trivial) methods
may be required to compute a factorization 
of the numerator of the resulting generating function.

We have found the techniques in this paper to be most useful when implemented
and
used in conjunction with a computer algebra package like Maple.
Through computer experiments,
we have been able to discover new identities, whose proof is beyond the
scope of these methods.
However, we show in a forthcoming article \cite{cs2002}
that methods inspired from Bousquet-M\'elou and Eriksson \cite{BME2}
and Aa Ya Yee \cite{yee1,yee2} 
can give straightforward proofs of these types of results.

\begin{small}
\bibliography{bib9}

\begin{thebibliography}{10}

\bibitem{Om1}
George~E. Andrews.
\newblock Mac{M}ahon's partition analysis. {I}. {T}he lecture hall partition
  theorem.
\newblock In {\em Mathematical essays in honor of Gian-Carlo Rota (Cambridge,
  MA, 1996)}, pages 1--22. Birkh\"auser Boston, Boston, MA, 1998.

\bibitem{Om2}
George~E. Andrews.
\newblock {M}ac{M}ahon's {P}artition {A}nalysis {II}: {F}undamental theorems.
\newblock {\em Annals of Combinatorics}, 4(3-4):327--338, 2000.

\bibitem{Om4}
George~E. Andrews and Peter Paule.
\newblock Mac{M}ahon's partition analysis. {I}{V}. {H}ypergeometric multisums.
\newblock {\em S\'em. Lothar. Combin.}, 42, 1999.
\newblock The Andrews Festschrift (Maratea, 1998).

\bibitem{Om6}
George~E. Andrews, Peter Paule, and Axel Riese.
\newblock {M}ac{M}ahon's {P}artition {A}nalysis {VI}: {A} new reduction
  theorem.
\newblock {\em submitted}, 2000.

\bibitem{Om8}
George~E. Andrews, Peter Paule, and Axel Riese.
\newblock {M}ac{M}ahon's {P}artition {A}nalysis {VIII}: {P}lane partitions
  diamonds.
\newblock {\em submitted}, 2000.

\bibitem{Om3}
George~E. Andrews, Peter Paule, and Axel Riese.
\newblock {M}ac{M}ahon's {P}artition {A}nalysis {III}: {T}he {O}mega package.
\newblock {\em European J, Comb.}, 22(7):887--904, 2001.

\bibitem{Om9}
George~E. Andrews, Peter Paule, and Axel Riese.
\newblock {M}ac{M}ahon's {P}artition {A}nalysis {IX}: $k$-gon partitions.
\newblock {\em Bull. Austral. Math. Soc.}, 64(2):321--329, 2001.

\bibitem{Om7}
George~E. Andrews, Peter Paule, and Axel Riese.
\newblock {M}ac{M}ahon's {P}artition {A}nalysis {VII}: {C}onstrained
  compositions.
\newblock {\em submitted}, 2001.

\bibitem{Om5}
George~E. Andrews, Peter Paule, Axel Riese, and Volker Strehl.
\newblock Mac{M}ahon's partition analysis {V}: {B}ijections, recursions, and
  magic squares.
\newblock In {\em Algebraic combinatorics and applications (G\"o\ss weinstein,
  1999)}, pages 1--39. Springer, Berlin, 2001.

\bibitem{BME1}
Mireille Bousquet-M\'elou and Kimmo Eriksson.
\newblock Lecture hall partitions.
\newblock {\em Ramanujan J.}, 1(1):101--111, 1997.

\bibitem{BME2}
Mireille Bousquet-M\'elou and Kimmo Eriksson.
\newblock Lecture hall partitions {II.}
\newblock {\em Ramanujan J.}, 1(2):165--185, 1997.

\bibitem{CCH}
Rod Canfield, Sylvie Corteel, and Pawel Hitczenko.
\newblock Random partitions with non negative $r^{th}$ differences.
\newblock {\em Adv. Applied Maths}, 27:298--317, 2001.

\bibitem{cs2002}
Sylvie Corteel and Carla~D. Savage.
\newblock {A}nti-lecture hall compositions.
\newblock {\em Discrete Math., to appear}, 2002.

\bibitem{hic}
D.~R. Hickerson.
\newblock A partition identity of the {E}uler type.
\newblock {\em Amer. Math. Monthly}, 81:627--629, 1974.

\bibitem{MaMo}
Percy~A. MacMahon.
\newblock {\em Combinatory analysis}.
\newblock Chelsea Publishing Co., New York, 1960.
\newblock Two volumes (bound as one).

\bibitem{St73}
Richard~P. Stanley.
\newblock Linear homogeneous {D}iophantine equations and magic labelings of
  graphs.
\newblock {\em Duke Math. J.}, 40:607--632, 1973.

\bibitem{yee1}
Ae~Ja Yee.
\newblock On combinatorics of lecture hall partitions.
\newblock {\em Ramanujan Journal}, 5:247--262, 2001.
\newblock Preprint.

\bibitem{yee2}
Ae~Ja Yee.
\newblock On the refined lecture hall theorem.
\newblock {\em Discrete Math.}, 248(1-3):293--298, 2002.

\bibitem{SC}
Doron Zeilberger.
\newblock Sylvie {C}orteel's one line proof of a partition theorem generated by
  {A}ndrews-{P}aule-{R}iese's computer.
\newblock {\em Shalosh B. Ekhad's and Doron Zeilberger's Very Own Journal},
  1998.

\end{thebibliography}
\bibliographystyle{plain}
\end{small}

\end{document}